\documentclass{article}
\usepackage{amsmath}
\usepackage{amssymb}
\usepackage{amscd}

\def\be{\begin{equation}}       \def\ee{\end{equation}}
\def\bd{\begin{displaymath}}    \def\ed{\end{displaymath}}
\def\beq{\begin{eqnarray}}      \def\eeq{\end{eqnarray}}
\def\bseq{\begin{eqnarray*}}    \def\eseq{\end{eqnarray*}}
\def\ba{\begin{array}}          \def\ea{\end{array}}
\def\ben{\begin{enumerate}}     \def\een{\end{enumerate}}

\def\lmt{\longmapsto}

\def\cop{\Delta}

\def\cnt{\varepsilon}
\def\ot{\otimes}

\def\GLqtwo{{GL}_{q}(2)}
\def\GLhtwo{{GL}_{h}(2)}

\def\GLpqtwo{{GL}_{p,q}(2)}
\def\GLhh'two{{GL}_{h,h'}(2)}

\def\Grs{{G}_{r,s}}
\def\Gmk{{G}_{m,k}}
\def\Grss'{{G}_{r}^{s,s'}}
\def\Gmkk'{{G}_{m}^{k,k'}}

\def\GLtwo{{GL}({2})}

\def\rinv{{r}^{-1}}
\def\sinv{{s}^{-1}}

\def\ident{{\bf 1}}

\begin{document}
\begin{center}
{\bf 
CONTRACTION OF THE $\Grs$ QUANTUM GROUP TO ITS NON STANDARD ANALOGUE
AND CORRESPONDING COLOURED QUANTUM GROUPS
}

\bigskip
\bigskip
\bigskip

{\large\bf Deepak Parashar and Roger J. McDermott
}

\bigskip
{\sl School of Computer and Mathematical Sciences\\
The Robert Gordon University\\
St. Andrew Street, Aberdeen AB25 1HG\\
United Kingdom}\\
{\sl 
email: deeps@scms.rgu.ac.uk, rm@scms.rgu.ac.uk
}

\bigskip
\bigskip
{\bf Abstract}
\end{center}
\medskip

The quantum group $G_{r,s}$ provides a realisation of the two parameter
quantum $\GLpqtwo$ which is known to be related to the two parameter non
standard $\GLhh'two$ group via the contraction method. We apply the
contraction procedure to $\Grs$ and obtain a new Jordanian quantum
group $\Gmk$.  Furthermore, we provide a realisation of $\GLhh'two$ in
terms of $\Gmk$. The contraction procedure is then extended to the
coloured quantum group $GL_{r}^{\lambda, \mu}(2)$ to yield a new
Jordanian quantum group $GL_{m}^{\lambda, \mu}(2)$. Both $\Grs$ and $\Gmk$
are then generalised to their coloured versions which inturn provide
similar realisations of $GL_{r}^{\lambda, \mu}(2)$ and $GL_{m}^{\lambda,
\mu}(2)$.
\par
{\bf PACS: 02.20.-a, 02.20.Qs}
\bigskip
\bigskip
\bigskip
\begin{center}
To appear in {\sl J. Math. Phys.}
\end{center}

\newpage
\section{Introduction}

In recent years, a lot of interest has been generated in the study of
Non Standard (or Jordanian) deformations of Lie groups and algebras. For
$GL(2)$, the Jordanian deformation (also known as the {\em h-deformation})
was initially introduced in [1,2] with its two parametric generalisation
given by Aghamohammadi in [3]. This was extended to the supersymmetric
case $GL(1/1)$ in [4]. At the algebra level, the non standard deformation
$U_{h}(sl(2))$ of $sl(2)$ was first proposed by Ohn [5], the universal
$R$-matrix was presented in [6-8] and irreducible representations studied
in [9,10]. A peculiar feature of this deformation is that the
corresponding $R$-matrix is triangular i.e. $R^{-1}=R_{21}$. The group
$GL(2)$ admits two distinct deformations with central determinant:
$\GLqtwo$ and $\GLhtwo$, and these are the only possible such deformations
(up to isomorphism) [11]. In [12], an observation was made that the
$h$-deformation could be obtained by a singular limit of a similarity
transformation from the $q$-deformations of the group $GL(2)$. Given this
contraction procedure, it would be useful to look for Jordanian
deformations of other $q$-groups. 
\par

The two parameter quantum group $\Grs$ was proposed by Basu-Mallick in
[13] as a particular quotient of the multiparameter $q$-deformation of
$GL(3)$. The structure of $\Grs$ is interesting because it contains the
one parameter $q$-deformation of $GL(2)$ as a Hopf subalgebra and also
gives a simple realisation of the quantum group $\GLpqtwo$ in terms of the
generators of $\Grs$. As an initial step in the further study of this
quantum group, the authors have recently shown [15] that the dual Hopf
algebra to $\Grs$ may be realised using the method described by Sudbery
[14] i.e. as the algebra of tangent vectors at the identity. As well as
this, a bicovariant differential calculus on $\Grs$ has also been
constructed by the authors [16]. In the present paper, we investigate the
contraction procedure on $\Grs$ in order to obtain its non standard
counterpart. The generators of the contracted structure are employed to
realise the two parameter non standard deformation $\GLhh'two$. This is
similar to what happens in the $q$-deformed case. Furthermore, we extend
the contraction procedure to the case of `coloured' quantum groups and
obtain a new single parameter non-standard quantum group 
$GL_{m}^{\lambda,\mu}(2)$. $\Gmk$ is then extended to its coloured
version, obtained by contraction from the coloured $\Grs$. These new
coloured extension also provide realisations of single parameter coloured
quantum and Jordanian deformation of $\GLtwo$.
\par

In Section II, the block diagonal form of the $R$-matrix of $\Grs$ is
presented. The contraction on this block diagonal $R$-matrix is
carried out in Section III yielding the new $R$-matrix for $\Gmk$. Section
IV defines a new non standard group $\Gmk$ and Section V provides a new
realisation of the well known $\GLhh'two$. Contraction for coloured
quantum groups and coloured extensions of $\Grs$ and $\Gmk$ are given in
Section VI with Section VII detailing the coloured realisations. In
Section VIII, we summarise our results and briefly discuss possible
further work. Appendices $A$,$B$ and $C$ summarise the Hopf algebra
relations for $\Grs$, $\GLhh'two$ and coloured $GL_{m}^{\lambda,\mu}(2)$
respectively.

\section{The $\Grs$ $R$-matrix}

The quantum group $\Grs$ is generated by the matrix of generators
\[
T=\begin{pmatrix}a&b&0\\c&d&0\\0&0&f\end{pmatrix}
\]

where the first four generators $a$, $b$, $c$ and $d$ form a Hopf
subalgebra which is isomorphic to $\GLqtwo$ quantum group with deformation
parameter $q=r^{-1}$. The two parameter $\GLpqtwo$ can also be realised
through the generators of this $\Grs$ Hopf algebra provided the sets of
deformation parameters $(p,q)$ and $(r,s)$ are related to each other in a
particular fashion. This quantum group can, therefore, be used to realise
both $\GLqtwo$ and $\GLpqtwo$ quantum groups. The expression for the
$R$-matrix of $\Grs$ is given in [13]. Explicitly,
this reads

\[
R(\Grs) = \begin{pmatrix}
r & 0 & 0 & 0 & 0 & 0 & 0 & 0 & 0\\
0 & 1 & 0 & r-r^{-1} & 0 & 0 & 0 & 0 & 0\\
0 & 0 & s & 0 & 0 & 0 & r-r^{-1} & 0 & 0\\
0 & 0 & 0 & 1 & 0 & 0 & 0 & 0 & 0\\
0 & 0 & 0 & 0 & r & 0 & 0 & 0 & 0\\
0 & 0 & 0 & 0 & 0 & 1 & 0 & r-r^{-1} & 0\\
0 & 0 & 0 & 0 & 0 & 0 & s^{-1} & 0 & 0\\
0 & 0 & 0 & 0 & 0 & 0 & 0 & 1 & 0\\
0 & 0 & 0 & 0 & 0 & 0 & 0 & 0 & r\
\end{pmatrix}
\]

with entries labelled in the usual numerical order $(11)$, $(12)$, $(13)$,
$(21)$, $(22)$, $(23)$, $(31)$, $(32)$, $(33)$. We start with the
observation, first made by Aschieri and Castellani [17,18 ], that if we
reorder the indices of this $R$-matrix with the elements in the order
$(11)$, $(12)$, $(21)$, $(22)$, $(13)$, $(23)$, $(31)$, $(32)$, $(33)$,
then we obtain a block matrix, say $R_{q}$ which is similar to the form of
the $\GLqtwo$ $R$-matrix with the $q$ in the $R^{11}_{11}$ position itself
replaced by the $\GLqtwo$ $R$-matrix.

\[
R_{q} = \begin{pmatrix}
R(GL_{r}(2)) & 0 & 0 & 0\\
0 & S & \lambda I & 0\\
0 & 0 & S^{-1} & 0\\
0 & 0 & 0 & r\
\end{pmatrix}
\]

where $R(GL_{r}(2))$ is the $4\times 4$ $R$-matrix for $\GLqtwo$ with
$q=r$, $\lambda=r-r^{-1}$, $I$ is the $2\times 2$ identity matrix and $S$
is the $2\times 2$ matrix 
$S = \left( \begin{smallmatrix}s&0\\0&1\end{smallmatrix} \right)$ 
where $r$ and $s$ are the deformation parameters. The zeroes are the zero
matrices of appropriate order. The usual block structure of the $R$-matrix
is clearly seen in this form.  It is straightforward to check that the
$RTT$- relations with this new $R$-matrix give the known $\Grs$
commutation relations.
\par

\section{$R$-matrix Contraction}

It is well known [12] that the non standard $R$-matrix $R_{h}(2)$
can be obtained from the $q$-deformed $R_{q}(2)$ as a singular limit of a
similarity transformation
\[
R_{h}(2) = \lim_{q\rightarrow 1}(g^{-1}\ot g^{-1}) R_{q}(2)(g\ot g)
\]
where $g=\left( \begin{smallmatrix}1&\eta\\0&1\end{smallmatrix} \right)$.
Such a transformation has been generalised to higher dimensions [19] and
has also been successfully applied to two parameter quantum groups.
Here we apply the above transformation for our $\Grs$ quantum group. Our
starting point is the block diagonal form of the $\Grs$ $R$-matrix,
denoted $R_{q}$
\[
R_{q} = \begin{pmatrix}
r & 0 & 0 & 0 & 0 & 0 & 0 & 0 & 0\\
0 & 1 & \lambda & 0 & 0 & 0 & 0 & 0 & 0\\
0 & 0 & 1 & 0 & 0 & 0 & 0 & 0 & 0\\
0 & 0 & 0 & r & 0 & 0 & 0 & 0 & 0\\
0 & 0 & 0 & 0 & s & 0 & \lambda & 0 & 0\\
0 & 0 & 0 & 0 & 0 & 1 & 0 & \lambda & 0\\
0 & 0 & 0 & 0 & 0 & 0 & s^{-1} & 0 & 0\\
0 & 0 & 0 & 0 & 0 & 0 & 0 & 1 & 0\\
0 & 0 & 0 & 0 & 0 & 0 & 0 & 0 & r
\end{pmatrix}
\]
where $\lambda=r-r^{-1}$. We apply to $R_{q}$ the transformation
\[
(G^{-1}\ot G^{-1})R_{q}(G\ot G)
\]
Here the transformation matrix $G$ is a $3\times 3$ matrix and chosen in
the block diagonal form
\[
G=\begin{pmatrix}g&0\\0&1\end{pmatrix}
\]
where $g$ is the transformation matrix for the two dimensional case. The
similarity transformation gives the matrix
\[
\begin{pmatrix}
r & \eta (r-1) & \eta (r-\lambda -1) & \eta^{2}(2(r-1)-\lambda) & 0 & 0 &
0 & 0 & 0\\
0 & 1 & \lambda & \eta (\lambda +1-r) & 0 & 0 & 0 & 0 & 0\\
0 & 0 & 1 & \eta (1-r) & 0 & 0 & 0 & 0 & 0\\
0 & 0 & 0 & r & 0 & 0 & 0 & 0 & 0\\
0 & 0 & 0 & 0 & s & \eta (s-1) & \lambda & 0 & 0\\
0 & 0 & 0 & 0 & 0 & 1 & 0 & \lambda & 0\\
0 & 0 & 0 & 0 & 0 & 0 & s^{-1} & \eta (s^{-1}-1) & 0\\
0 & 0 & 0 & 0 & 0 & 0 & 0 & 1 & 0\\
0 & 0 & 0 & 0 & 0 & 0 & 0 & 0 & r
\end{pmatrix}
\]
Following the procedure outlined in [12], we substitute
$\eta=\frac{m}{1-r}$ to obtain
\[
\begin{pmatrix}
r & -m & mr^{-1} & m^{2}r^{-1} & 0 & 0 & 0 & 0 & 0\\
0 & 1 & r-r^{-1} & -mr^{-1} & 0 & 0 & 0 & 0 & 0\\
0 & 0 & 1 & m & 0 & 0 & 0 & 0 & 0\\
0 & 0 & 0 & r & 0 & 0 & 0 & 0 & 0\\
0 & 0 & 0 & 0 & s & -m\frac{1-s}{1-r} & r-r^{-1} & 0 & 0\\
0 & 0 & 0 & 0 & 0 & 1 & 0 & r-r^{-1} & 0\\
0 & 0 & 0 & 0 & 0 & 0 & s^{-1} & ms^{-1}\frac{1-s}{1-r} & 0\\
0 & 0 & 0 & 0 & 0 & 0 & 0 & 1 & 0\\
0 & 0 & 0 & 0 & 0 & 0 & 0 & 0 & r
\end{pmatrix}
\]
In the limit $r\rightarrow 1$, $s\rightarrow 1$ such that
$\frac{1-s}{1-r}\rightarrow \frac{k}{m}$, this yields the Jordanian 
$R$-matrix
\[
R_{h} = R(\Gmk) = \begin{pmatrix}
1 & -m & m & m^{2} & 0 & 0 & 0 & 0 & 0\\
0 & 1 & 0 & -m & 0 & 0 & 0 & 0 & 0\\
0 & 0 & 1 & m & 0 & 0 & 0 & 0 & 0\\
0 & 0 & 0 & 1 & 0 & 0 & 0 & 0 & 0\\
0 & 0 & 0 & 0 & 1 & k & 0 & 0 & 0\\
0 & 0 & 0 & 0 & 0 & 1 & 0 & 0 & 0\\
0 & 0 & 0 & 0 & 0 & 0 & 1 & -k & 0\\
0 & 0 & 0 & 0 & 0 & 0 & 0 & 1 & 0\\
0 & 0 & 0 & 0 & 0 & 0 & 0 & 0 & 1
\end{pmatrix}
\]
where the entries are labelled in the block diagonal form $(11)$, $(12)$,
$(21)$, $(22)$, $(13)$, $(23)$, $(31)$, $(32)$, $(33)$. It is
straightforward to verify that this $R$-matrix is triangular and a
solution of the Quantum Yang Baxter Equation
\[
R_{12}R_{13}R_{23} = R_{23}R_{13}R_{12}
\]
It is interesting to note that the block diagonal form of $R(\Gmk)$ embeds
in the top left corner the $R$-matrix for the single parameter deformed
$\GLhtwo$ for $m=h$.

\section{The Non Standard $\Gmk$}

The contracted $R$-matrix $R(\Gmk)$ can be used in conjunction with
a $T$-matrix of generators of the form
\[
T=\begin{pmatrix}a&b&0\\c&d&0\\0&0&f\end{pmatrix}
\]
to form a two parameter non standard quantum group $\Gmk$. The
$RTT$- relations
\[
RT_{1}T_{2} = T_{2}T_{1}R
\]
(where $T_{1}=T\ot \ident$, $T_{2}=\ident\ot T$) give the commutation
relations between the generators $a$, $b$, $c$, $d$ and $f$.
\[
[c,d] = -mc^{2}, \qquad[c,b] = -m(ac+cd) = -m(ca+dc)
\]
\[
[c,a] = -mc^{2}, \qquad[d,a] = -m(d-a)c = -mc(d-a)
\]
\[
[d,b] = -m(d^{2}-\delta)
\]
\[
[b,a] = -m(\delta-a^{2})
\]
and
\[
[f,a] = kcf, \qquad[f,b] = k(df-fa)
\]
\[
[f,c] = 0, \qquad[f,d] = -kcf
\]
The element $\delta = ad-bc+mac = ad-cb-mcd$ is central in the whole
algebra.
(Note that the first set of the above relations consist of
elements $a$, $b$, $c$ and $d$ which form a subalgebra that coincides
exactly with the single parameter non standard $\GLhtwo$ for $m=h$.) The
coalgebra structure of $\Gmk$ can be written as
\[
\begin{array}{l}
\cop\begin{pmatrix}a&b&0\\c&d&0\\0&0&f\end{pmatrix} = 
\begin{pmatrix}
a\ot a+b\ot c & a\ot b+b\ot d & 0\\
c\ot a+d\ot c & c\ot b+d\ot d & 0\\
0 & 0 & f\ot f
\end{pmatrix}\\
\\
\cnt\begin{pmatrix}a&b&0\\c&d&0\\0&0&f\end{pmatrix} = 
\begin{pmatrix}1&0&0\\0&1&0\\0&0&1\end{pmatrix}
\end{array}
\]

Adjoining the element $\delta^{-1}$ to the algebra enables determination
of the antipode matrix $S(T)$,
\[
S\begin{pmatrix}a&b&0\\c&d&0\\0&0&f\end{pmatrix}
=\delta^{-1}
\begin{pmatrix}
d+mc & -b+m(d-a)+m^{2}c & 0\\
-c & a-mc & 0\\
0&0&\delta f^{-1}\end{pmatrix}
\]

(The Hopf structure of $\delta^{-1}$ is $\cop(\delta^{-1})=\delta^{-1}\ot
\delta^{-1}, \cnt(\delta^{-1})=\ident, S(\delta^{-1})=\delta$.) These
relations are consistent with the usual axioms of the Hopf algebra,
\begin{eqnarray*}
\mu\circ (id\ot \mu)=\mu\circ (\mu\ot id)&,& (id\ot \cop)\circ
\cop=(\cop\ot id)\circ \cop\\
(id\ot \cnt)\circ \cop=(\cnt\ot id)\circ \cop=id&,& \mu\circ (id\ot
S)\circ \cop=\mu\circ (S\ot
id)\circ \cop=\ident\circ \cnt\\
\cop(xy)=\cop(x)\cop(y), \qquad \cnt(xy)&=&\cnt(x)\cnt(y), \qquad
S(xy)=S(y)S(x)
\end{eqnarray*}

where $\mu$ denotes the multiplication operation $\mu (x\ot y)=xy$ and
$id$ is the identity transformation. 
\par

It is evident that the elements $a$, $b$, $c$ and $d$ of $\Gmk$ form a
Hopf subalgebra which coincides with non standard $GL(2)$ with deformation
parameter $m$. This is exactly analogous to the $q$-deformed case where
the first four elements of $\Grs$ form the $\GLqtwo$ Hopf subalgebra.
Again, the remaining fifth element $f$ generates the $GL(1)$ group, as it
did in the $q$-deformed case, and the second parameter appears only
through the cross commutation relations between $GL_{m}(2)$ and $GL(1)$
elements. Therefore, $\Gmk$ can also be considered as a two parameter
Jordanian deformation of classical $GL(2)\ot GL(1)$ group.
\par

\section{A Realisation of $\GLhh'two$}

Now we investigate the connection of the newly defined $\Gmk$
with the non standard two parameter $\GLhh'two$. It was observed by
Basu-Mallick [13] that there is a Hopf algebra homomorphism $\mathcal{F}$
from $\Grs$ to $\GLpqtwo$ given by
\[
\mathcal{F}_{N}: \Grs\lmt \GLpqtwo
\]
\[
\mathcal{F}_{N}:\begin{pmatrix}a&b\\c&d\end{pmatrix}\lmt
\begin{pmatrix}a'&b'\\c'&d'\end{pmatrix}=
f^{N}\begin{pmatrix}a&b\\c&d\end{pmatrix}
\]

The elements $a'$,$b'$,$c'$ and $d'$ are the generators of $\GLpqtwo$ and
$N$ is a fixed non-zero integer. The relation between the deformation
parameters $(p,q)$ and $(r,s)$ is given by
\[p = \rinv s^{N} \quad \text{and} \quad q = \rinv s^{-N}\]

A Hopf algebra homomorphism of exactly the same form exists between the
generators of $\Gmk$ and $\GLhh'two$ which is straightforward to verify.
Moreover, the two sets of deformation parameters $(h,h')$ and $(m,k)$ are
related via the equation
\[
m = -h+Nk = -h'-Nk
\]
\[
\text{i.e.} \quad \quad h = -m + Nk \quad \text{and} \quad h' = -m - Nk
\]
Note that for vanishing $k$, have $h'=h$ and one gets the one
parameter case. In addition, using the above realisation together with the
coproduct, counit and antipode axioms for the $\Gmk$ algebra and the
respective homeomorphism properties, one can easily recover the standard
coproduct, counit and antipode for $\GLhh'two$. Thus, the non
standard $\GLhh'two$ group can in fact be reproduced from the newly
defined non standard $\Gmk$. The above realisations can be exhibited in
the following commutative diagram
\[
\newcommand{\contr}{\operatorname{contraction}}
\begin{CD}
\Grs @>\mathcal{F}>> \GLpqtwo\\
@V{\contr}VV	@VV{\contr}V\\
\Gmk @>>\mathcal{F}> \GLhh'two
\end{CD}
\]

It is curious to note that if we write $p=e^{h}$, $q=e^{h'}$, $r=e^{m}$
and $s=e^{k}$, then the relations between the parameters in the
$q$-deformed case and the $h$-deformed case are identical. 
\par

\section{Coloured Quantum Groups}

The standard quantum group relations can be extended by parametrising the
corresponding generators using some continuous `colour' variables and
redefining the associated algebra and coalgebra in a way that all Hopf
algebraic properties remain preserved [13, 20, 21]. For the case of a
single parameter quantum deformation of $\GLtwo$ (with deformation
parameter $r$), its `coloured' version [13] is given by the $R$-matrix

\[
R_{r}^{\lambda,\mu} = \begin{pmatrix}
r^{1-(\lambda-\mu)} & 0 & 0 & 0\\
0 & r^{\lambda+\mu} & r-r^{-1} & 0\\
0 & 0 & r^{-(\lambda+\mu)} & 0\\
0 & 0 & 0 & r^{1+(\lambda-\mu)}\
\end{pmatrix}
\]
which satisfies
\[
R_{12}^{\lambda,\mu}R_{13}^{\lambda,\nu}R_{23}^{\mu,\nu} =
R_{23}^{\mu,\nu}R_{13}^{\lambda,\nu}R_{12}^{\lambda,\mu}
\]
the so-called `Coloured' Quantum Yang Baxter Equation (CQYBE). This gives
rise to the coloured $RTT$ relations
\[
R_{r}^{\lambda,\mu}T_{1\lambda}T_{2\mu}=T_{2\mu}T_{1\lambda}
R_{r}^{\lambda,\mu}
\]
(where $T_{1\lambda}=T_{\lambda}\ot \ident$ and $T_{2\mu}=\ident\ot
T_{\mu}$) in which the entries of the $T$ matrices carry colour 
dependence. The coproduct and counit for the coalgebra structure are given
by
\begin{eqnarray*}
\cop (T_{\lambda})&=&T_{\lambda}\ot T_{\lambda}\\
\cnt (T_{\lambda})&=&\ident
\end{eqnarray*}
and depend only on one colour parameter. By contrast, the algebra 
structure is more complicated with only generators of two different
colours appearing simultaneously in the algebraic relations. The full Hopf
algebraic structure can be constructed and results in a coloured extension
of the quantum group. Since $\lambda$ and $\mu$ are continuous variables,
this implies the coloured quantum group has an infinite number of
generators. The quantum determinant
$D_{\lambda}=a_{\lambda}d_{\lambda}-r^{-(1+2\lambda)}c_{\lambda}b_{\lambda}$
is group-like but not central, and the antipode is
\[
S(T_{\lambda})= 
D_{\lambda}^{-1}\begin{pmatrix}d_{\lambda}&-r^{1+2\lambda}b_{\lambda}\\
-r^{-1-2\lambda}&a_{\lambda}\end{pmatrix}
\]

In order to investigate the contraction for coloured quantum groups. we
apply to $R_{r}^{\lambda,\mu}$ the transformation
\[
(g\ot g)^{-1}R_{r}^{\lambda,\mu}(g\ot g)
\]
where $g$ is the two dimensional transformation matrix
$\left( \begin{smallmatrix}1&\eta\\0&1\end{smallmatrix} \right)$ and
$\eta$ is chosen to be $\eta=\frac{m}{1-r}$. This gives the matrix
\[
\begin{pmatrix}
r^{1-(\lambda-\mu)} & R_{12}^{11} &
R_{21}^{11} & R_{22}^{11}\\
0 & r^{\lambda+\mu} & r-r^{-1} & R_{22}^{12}\\   
0 & 0 & r^{-(\lambda+\mu)} & R_{22}^{21}\\
0 & 0 & 0 & r^{1+(\lambda-\mu)}
\end{pmatrix}
\]
where
\[
\begin{array}{ll}
&R_{12}^{11}=mr^{1-\lambda +\mu}[-1+2\lambda]_{r}\\
&R_{21}^{11}=-mr^{1-\lambda +\mu}[1+2\mu]_{r}+mr^{-1}(1+r)\\
&R_{22}^{12}=mr^{\lambda+ \mu}[1-2\mu]_{r}-mr^{-1}(1+r)\\
&R_{22}^{21}=mr^{-(\lambda +\mu)}[1+2\lambda]_{r}\\
&R_{22}^{11}= 
-m^{2}r([-\lambda+\mu]_{r}[-1+2\lambda]_{r}+[\lambda-\mu]_{r}[-1-2\lambda]_{r})
+ m^{2}r^{-1}([1+2\lambda]_{r}[1-2\lambda]_{r})
\end{array}
\]
and $[x]_{r}=(\frac{1-x^{r}}{1-x})$ denotes the basic number from
$q$-analysis.

In the limit $r\rightarrow 1$, we obtain a new $R$-matrix
\[
R_{m}^{\lambda,\mu} = \begin{pmatrix}
1 & -m(1-2\lambda) & m(1-2\mu) & m^{2}(1-4\lambda \mu)\\
0 & 1 & 0 & -m(1+2\mu)\\
0 & 0 & 1 & m(1+2\lambda)\\
0 & 0 & 0 & 1
\end{pmatrix}
\]
which is a coloured $R$-matrix for a Jordanian deformation of $\GLtwo$.
This $R$-matrix satisfies the CQYBE and is `colour' triangular i.e.
$R_{12}^{\lambda,\mu}=({R_{21}^{\mu,\lambda}})^{-1}$, a coloured extension
of the notion of triangularity. This $R$-matrix is distinct from that of
the coloured Jordanian deformation of $\GLqtwo$ obtained in [20,21] by
other means. The Hopf algebra structure and the commutation relations for
the quantum group associated with this $R$-matrix are given in Appendix
$C$.
\par
{\bf Note:} This is the first time that such a contraction procedure has
been
applied to obtain a coloured Jordanian $R$-matrix and hence the coloured
Jordanian quantum group.

\subsection{Coloured Extension of $\Grs$ : $\Grss'$}

The coloured extension of $\Grs$ proposed in [13] has only one deformation
parameter $r$ and two colour paramters $s$ and $s'$. The second
deformation parameter of the uncoloured case now plays the role of a
colour parameter. In such a coloured extension, the first four
generators $a$,$b$,$c$,$d$ are kept independent of the colour parameter(s)
while the fifth generator $f$ is now paramterised by $s$ and $s'$.
The matrices of generators are
\[
T_{s}=\begin{pmatrix}a&b&0\\c&d&0\\0&0&f_{s}\end{pmatrix}\quad , \quad
T_{s'}=\begin{pmatrix}a&b&0\\c&d&0\\0&0&f_{s'}\end{pmatrix}
\]
From the $RTT$ relations, one observes that the commutation relations
between $a$,$b$,$c$,$d$ are as before but $f_{s}$ and $f_{s'}$ now satisfy
\[
\begin{array}{llll}
af_{s}=f_{s}a,\quad &bf_{s}=s^{-1}f_{s}b,\quad &cf_{s}=sf_{s}c,\quad 
&df_{s}=f_{s}d\\
af_{s'}=f_{s'}a,\quad &bf_{s'}={s'}^{-1}f_{s'}b,\quad 
&cf_{s'}=s'f_{s'}c,\quad &df_{s'}=f_{s'}d
\end{array}
\]
and
\[
[f_{s},f_{s'}]=0
\]
Our choice of the $9\times 9$ $R$-matrix for $\Grss'$ is

\[
R_{r}^{s,s'} = \begin{pmatrix}
R_{r}(2) & 0 & 0 & 0\\
0 & S(s,s') & 0 & 0\\
0 & 0 & \overline{S}(s,s') & 0\\
0 & 0 & 0 & r
\end{pmatrix}
\]
where
\[
R_{r}(2)=\begin{pmatrix}
r & 0 & 0 & 0\\                    
0 & 1 & r-r^{-1} & 0\\         
0 & 0 & 1 & 0\\            
0 & 0 & 0 & r
\end{pmatrix}, \quad
S(s,s')=\begin{pmatrix}
\sqrt{ss'} & 0\\
0 & \sqrt{\frac{s}{s'}}
\end{pmatrix}, \quad
\overline{S}(s,s')=\begin{pmatrix}
\frac{1}{\sqrt{ss'}} & 0\\
0 & \sqrt{\frac{s}{s'}}
\end{pmatrix}
\]
which satisfies the CQYBE
\[
R_{12}(r;s,s')R_{13}(r;s,s'')R_{23}(r;s',s'') =
R_{23}(r;s',s'')R_{13}(r;s,s'')R_{12}(r;s,s')
\]
$S$ and $\overline{S}$ satisfy the exchange relation
$\overline{S}(s,s')={S(s',s)}^{-1}$.

\subsection{Coloured Extension of $\Gmk$ : $\Gmkk'$}

Similar to the case of $\Grs$, we propose a coloured extension of the
Jordanian quantum group $\Gmk$. The first four generators remain
independent of the coloured parameters $k$ and $k'$ whereas the generator
$f$ is parameterised by $k$ and $k'$. Again, the second
deformation parameter $k$ of the uncoloured case now plays the role of a
colour parameter and the $T$-matrices are
\[
T_{k}=\begin{pmatrix}a&b&0\\c&d&0\\0&0&f_{k}\end{pmatrix}\quad , \quad
T_{k'}=\begin{pmatrix}a&b&0\\c&d&0\\0&0&f_{k'}\end{pmatrix}
\]
The commutation relations between $a$,$b$,$c$,$d$ remain unchanged whereas
$f_{k}$ and $f_{k'}$ satisfy
\[
\begin{array}{lllll}
&[f_{k},a]=kcf_{k},\quad &[f_{k},b]=k(df_{k}-f_{k}a),\quad
&[f_{k},c]=0,\quad &[f_{k},d]=-kcf_{k}\\
&[f_{k'},a]=k'cf_{k'},\quad &[f_{k'},b]=k'(df_{k'}-f_{k'}a),\quad
&[f_{k'},c]=0,\quad &[f_{k'},d]=-k'cf_{k'}
\end{array}
\]
and
\[
[f_{k},f_{k'}]=0
\]
Our choice of the $9\times 9$ $R$-matrix for $\Gmkk'$ is

\[
R_{m}^{k,k'} = \begin{pmatrix}
R_{m}(2) & 0 & 0 & 0\\
0 & K(k,k') & 0 & 0\\
0 & 0 & \overline{K}(k,k') & 0\\
0 & 0 & 0 & 1
\end{pmatrix}
\]
where
\[
R_{m}(2)=\begin{pmatrix}
1 & -m & m & m^{2}\\                    
0 & 1 & 0 & -m\\         
0 & 0 & 1 & m\\            
0 & 0 & 0 & 1
\end{pmatrix}, \quad
K(k,k')=\begin{pmatrix}
1 & -k'\\
0 & 1
\end{pmatrix}, \quad
\overline{K}(k,k')=\begin{pmatrix}
1 & k\\
0 & 1
\end{pmatrix}
\]
This $R$-matrix is chosen since it is the contraction limit of the 
$R$-matrix for the coloured extension of $\Grs$ via the transformation
\[
R_{m}^{k,k'} = \lim_{r\rightarrow 1}(G\ot G)^{-1} R_{r}^{s,s'}(G\ot G)
\]
where
\[
G=\begin{pmatrix}g&0\\0&1\end{pmatrix};\quad
g=\begin{pmatrix}1&\eta\\0&1\end{pmatrix},\quad \eta=\frac{m}{r-1}
\]
It is a solution of the CQYBE
\[
R_{12}(m;k,k')R_{13}(m;k,k'')R_{23}(m;k',k'') =
R_{23}(m;k',k'')R_{13}(m;k,k'')R_{12}(m;k,k')
\]
and is colour triangular. Again, $K$ and $\overline{K}$ satisfy the
exchange relation $\overline{K}(k,k')={K(k',k)}^{-1}$. 

\section{Coloured Realisations}

It is shown in [13] that, similar to the uncoloured case, the coloured
quantum group $\Grss'$ provides a realisation of the well known coloured
$GL_{r}^{\lambda,\mu}(2)$ where Hopf algebra homomorphism from $\Grss'$
to $GL_{r}^{\lambda,\mu}(2)$
\[
\mathcal{F}_{N}: \Grss'\lmt GL_{r}^{\lambda,\mu}(2)
\]
is given by
\[
\mathcal{F}_{N}:\begin{pmatrix}a&b\\c&d\end{pmatrix}\lmt
\begin{pmatrix}a'_{\lambda}&b'_{\lambda}\\c'_{\lambda}&d'_{\lambda}
\end{pmatrix}=f_{s}^{N}\begin{pmatrix}a&b\\c&d\end{pmatrix}
\]
\[
\mathcal{F}_{N}:\begin{pmatrix}a&b\\c&d\end{pmatrix}\lmt
\begin{pmatrix}a'_{\mu}&b'_{\mu}\\c'_{\mu}&d'_{\mu}
\end{pmatrix}=f_{s'}^{N}\begin{pmatrix}a&b\\c&d\end{pmatrix}
\]
where N is a fixed non-zero integer and the sets of colour parameters
$(s,s')$ and $(\lambda,\mu)$ are related through quantum deformation
parameter $r$ by
\[
s=r^{2N\lambda} \quad , \quad s'=r^{2N\mu}
\]
The primed generators $a'_{\lambda}$,
$b'_{\lambda}$,$c'_{\lambda}$,$d'_{\lambda}$ and
$a'_{\mu}$,$b'_{\mu}$,$c'_{\mu}$,$d'_{\mu}$ belong to
$GL_{r}^{\lambda,\mu}(2)$ whereas the unprimed ones
$a$,$b$,$c$,$d$,$f_{s}$ and  $f_{s'}$ are generators of $\Grss'$.
\par
Just as in the $q$-deformed case, we obtain a realisation of
the $h$-deformed quantum group $GL_{m}^{\lambda,\mu}(2)$ using the newly
defined coloured quantum group $\Gmkk'$. If we again denote the generators
of   
$GL_{m}^{\lambda,\mu}(2)$ by
$a'_{\lambda}$,$b'_{\lambda}$,$c'_{\lambda}$,$d'_{\lambda}$ and 
$a'_{\mu}$,$b'_{\mu}$,$c'_{\mu}$,$d'_{\mu}$ and the generators of $\Gmkk'$
by $a$,$b$,$c$,$d$,$f_{k}$ and  $f_{k'}$ then a Hopf algebra homomorphism
from $\Gmkk'$ to
$GL_{m}^{\lambda,\mu}(2)$
\[
\mathcal{F}_{N}: \Gmkk'\lmt GL_{m}^{\lambda,\mu}(2)
\]
is of exactly the same form
\[
\mathcal{F}_{N}:\begin{pmatrix}a&b\\c&d\end{pmatrix}\lmt
\begin{pmatrix}a'_{\lambda}&b'_{\lambda}\\c'_{\lambda}&d'_{\lambda}
\end{pmatrix}=f_{k}^{N}\begin{pmatrix}a&b\\c&d\end{pmatrix}
\]   
\[
\mathcal{F}_{N}:\begin{pmatrix}a&b\\c&d\end{pmatrix}\lmt
\begin{pmatrix}a'_{\mu}&b'_{\mu}\\c'_{\mu}&d'_{\mu}
\end{pmatrix}=f_{k'}^{N}\begin{pmatrix}a&b\\c&d\end{pmatrix}
\]
The sets of colour parameters $(k,k')$ and $(\lambda,\mu)$ are related
to the Jordanian deformation parameter $m$ by
\[
Nk=-2m\lambda \quad , \quad Nk'=-2m\mu
\]   
and N, again, is a fixed non-zero integer. 
\par

\section{Conclusions}

In this work, we have applied the contraction procedure to the $\Grs$
quantum group and obtained a new Jordanian quantum group $\Gmk$. The
group $\Gmk$ has five generators and two deformation parameters and
contains the single parameter $\GLhtwo$ as a Hopf subalgebra (generated by
the first four elements). The remaining fifth generator corresponds to the
$GL(1)$ group. Furthermore, we have given a realisation of the two
parameter $\GLhh'two$ through the generators of $\Gmk$ which
also reproduces its full Hopf algebra structure. The results match with
the $q$-deformed case. The bigger picture that emerges from our analysis
of contraction for uncoloured as well as coloured quantum groups and their
morphisms can be represented in the following diagram
\[
\begin{CD}
\GLpqtwo @<\mathcal{F}<< \Grs @>\mathcal{E}>> \Grss' @>\mathcal{F}>>
GL_{r}^{\lambda,\mu}(2)\\
@V{\mathcal{C}}VV    @VV{\mathcal{C}}V @VV{\mathcal{C}}V
@VV{\mathcal{C}}V\\
\GLhh'two @<<\mathcal{F}< \Gmk @>>\mathcal{E}> \Gmkk' @>>\mathcal{F}>
GL_{m}^{\lambda,\mu}(2)
\end{CD}
\]
where $\mathcal{C}$, $\mathcal{F}$ and $\mathcal{E}$ denote the
contraction, Hopf algebra homomorphism and coloured extension 
respectively. The objects at the top level are the $q$ deformed ones and
the corresponding Jordanian counterparts are shown at the bottom level
of the diagram.
\par
Future work along the lines indicated in the paper will lead to
the explicit derivation of the dual algebra for the non standard $\Gmk$
and its coloured extension, and this will be presented by the authors in
a later paper.

\newpage

\section*{References}

[1] E. E. Demidov, Yu. I. Manin, E. E. Mukhin and D. V. Zhdanovich, Prog.
Theor. Phys. Suppl. {\bf 102}, 203 (1990).\par
[2] S. Zakrzewski, Lett. Math. Phys. {\bf 22}, 287 (1991).\par
[3] A. Aghamohammadi, Mod. Phys. Lett. {\bf A8}, 2607 (1993).\par
[4] L. Dabrowski and P. Parashar, Lett. Math. Phys. {\bf 38}, 331
(1996).\par
[5] C. Ohn, Lett. Math. Phys. {\bf 25}, 85 (1992).\par
[6] A. A. Vladimirov, Phys. Lett. {\bf A8}, 2573 (1993).\par
[7] A. Shariati, A. Aghamohammadi and M. Khorrami, Mod. Phys. Lett. {\bf
A11}, 187 (1996).\par
[8] A. Ballesteros and F. J. Herranz, J. Phys. {\bf A29}, L311 (1996).\par
[9] V. K. Dobrev, ICTP Preprint - IC/96/14 (1996).\par
[10] B. Abdessalam, A. Chakrabarti and R. Chakrabarti, Mod. Phys. Lett.
{\bf A11}, 2883 (1996).\par
[11] B. A. Kupershmidt, J. Phys. {\bf A25}, L1239 (1992).\par
[12] A. Aghamohammadi, M. Khorrami and A. Shariati, J. Phys. {\bf A28},
L225 (1995).\par
[13] B. Basu-Mallick, hep-th/9402142; Intl. J. Mod. Phys. {\bf A10}, 2851
(1995).\par
[14] A. Sudbery, Proc. Workshop on Quantum Groups, Argonne (1990)
eds. T. Curtright, D. Fairlie   and C. Zachos, pp. 33-51.\par
[15] D. Parashar and R. J. McDermott, math.QA/9901132.\par
[16] D. Parashar and R. J. McDermott, RGU Preprint (1999).\par
[17] P. Aschieri and L. Castellani, Intl. J. Mod. Phys. {\bf A11},
1019 (1996).\par
[18] P. Aschieri, math.QA/9805119.\par
[19] M. Alishahiha, J. Phys. {\bf A28}, 6187 (1995).\par
[20] P. Parashar, Lett. Math. Phys. {\bf 45}, 105 (1998).\par
[21] C. Quesne, J. Math. Phys. {\bf 38}, 6018 (1997); {\em ibid} {\bf 
39}, 1199 (1998).\par

\newpage

\section*{Appendix $A$ - $\Grs$ Quantum Group}

The two parameter quantum group $\Grs$ is generated by elements $a$, $b$,
$c$, $d$, and
$f$ satisfying the relations
\[
\begin{array}{ll}
ab=\rinv ba,&db=rbd\\
ac=\rinv ca,&dc=rcd\\
bc=cb,&[a,d]=(\rinv-r)bc
\end{array}
\]
and
\[
\begin{array}{ll}
af=fa,&cf=sfc\\
bf=\sinv fb,&df=fd
\end{array}
\]
Elements $a$, $b$, $c$, $d$ satisfying the first set of
commutation relations form a subalgebra which coincides exactly with
$\GLqtwo$ when $q = \rinv$.The Hopf structure is given as
\[
\begin{array}{l}
\cop\begin{pmatrix}a&b&0\\c&d&0\\0&0&f\end{pmatrix} = 
\begin{pmatrix}
a\ot a+b\ot c & a\ot b+b\ot d & 0\\
c\ot a+d\ot c & c\ot b+d\ot d & 0\\
0 & 0 & f\ot f
\end{pmatrix}\\
\\
\cnt\begin{pmatrix}a&b&0\\c&d&0\\0&0&f\end{pmatrix} = 
\begin{pmatrix}1&0&0\\0&1&0\\0&0&1\end{pmatrix}
\end{array}
\]
The Casimir operator is defined as $\delta = ad-\rinv bc$. The inverse is
assumed to exist and satisfies $\cop(\delta^{-1})=\delta^{-1}\ot 
\delta^{-1}$, $\cnt (\delta^{-1})=1$, $S(\delta^{-1})=\delta$, which
enables determination of the antipode matrix $S(T)$, as
\[
S\begin{pmatrix}a&b&0\\c&d&0\\0&0&f\end{pmatrix}
=\delta^{-1}
\begin{pmatrix}d&-rc&0\\-\rinv c&a&0\\0&0&\delta f^{-1}\end{pmatrix}
\]
The quantum determinant $D = \delta f$ is group-like but not central.

\newpage
\section*{Appendix $B$ - Non Standard $\GLhh'two$}

The two parameter non standard group $\GLhh'two$ is generated by the
matrix of generators $T=\left(
\begin{smallmatrix}a&b\\c&d\end{smallmatrix} \right)$ and the $4\times 4$
$R$-matrix is given as
\[
R = \begin{pmatrix}
1 & -h' & h' & hh'\\
0 & 1 & 0 & -h\\
0 & 0 & 1 & h\\
0 & 0 & 0 & 1
\end{pmatrix}
\]
The commutation relations among the generators $a$, $b$, $c$ and $d$ are
\[
[a,c] = hc^{2}, \qquad[b,c] = hcd+h'ac
\]
\[
[d,c] = h'c^{2}, \qquad[a,d] = hcd-h'ca
\]
\[
[a,b] = h'(D-a^{2})
\]
\[
[d,b] = h(D-d^{2})
\]
where $D = ad-cb-hcd = ad-bc+h'ac$ is the quantum determinant and $h$ and
$h'$ are the two deformation parameters. The Hopf algebra structure is
\[
\begin{array}{l}
\cop\begin{pmatrix}a&b\\c&d\end{pmatrix} = 
\begin{pmatrix}
a\ot a+b\ot c & a\ot b+b\ot d\\
c\ot a+d\ot c & c\ot b+d\ot d
\end{pmatrix}\\
\\
\cnt\begin{pmatrix}a&b\\c&d\end{pmatrix} = 
\begin{pmatrix}1&0\\0&1\end{pmatrix}
\end{array}
\]
$D^{-1}$ exists and satisfies $\cop(D^{-1})=D^{-1}\ot D^{-1},
\cnt(D^{-1})=\ident, S(D^{-1})=D$. Using this structure the antipode
matrix can be expressed as
\[
S(T) = D^{-1} \begin{pmatrix}
d+hc & -b+h(d-a)+h^{2}c\\
-c & a-hc\\
\end{pmatrix}
=\begin{pmatrix} 
d+h'c & -b+h'(d-a)+{h'}^{2}c\\    
-c & a-h'c\\    
\end{pmatrix} D^{-1}
\]

\newpage

\section*{Appendix $C$ - Coloured Jordanian $GL_{m}^{\lambda,\mu}(2)$}

The commutation relations between the generating elements
$a_{\lambda}$,$b_{\lambda}$,$c_{\lambda}$,$d_{\lambda}$ and
$a_{\mu}$,$b_{\mu}$,$c_{\mu}$,$d_{\mu}$ of
$GL_{m}^{\lambda,\mu}(2)$ are
\[
\begin{array}{ll}
&[a_{\lambda},c_{\mu}]= m(1+2\mu)c_{\lambda}c_{\mu}\\
&[d_{\lambda},c_{\mu}]= m(1-2\mu)c_{\mu}c_{\lambda}\\
&[a_{\lambda},d_{\mu}]=
m(1+2\mu)c_{\lambda}d_{\mu}-m(1-2\lambda)c_{\mu}a_{\lambda}\\
&[a_{\lambda},b_{\mu}]=
m(1-2\lambda)a_{\lambda}d_{\mu}-m(1-2\lambda)a_{\mu}a_{\lambda}-
m(1-2\mu)c_{\lambda}b_{\mu}-m^{2}(1-4\lambda \mu)c_{\lambda}d_{\mu}\\
&[b_{\lambda},d_{\mu}]=
m(1+2\mu)d_{\lambda}d_{\mu}+m(1+2\lambda)c_{\mu}b_{\lambda}-
m(1+2\mu)d_{\mu}a_{\lambda}+m^{2}(1-4\lambda \mu)c_{\mu}a_{\lambda}\\
&[b_{\lambda},c_{\mu}]=
m(1+2\mu)d_{\lambda}c_{\mu}+m(1-2\mu)c_{\mu}a_{\lambda}
\end{array}
\]
and
\[
\begin{array}{llll}
&[a_{\lambda},a_{\mu}]&=
&m(1-2\lambda)a_{\lambda}c_{\mu}-
m(1-2\mu)c_{\lambda}a_{\mu}-m^{2}(1-4\lambda \mu)c_{\lambda}c_{\mu}\\
&[b_{\lambda},b_{\mu}]&=
&m(1-2\lambda)b_{\lambda}d_{\mu}+m(1+2\lambda)a_{\mu}b_{\lambda}-
m(1-2\mu)d_{\lambda}b_{\mu}-m(1+2\mu)b_{\mu}a_{\lambda}\\
& & &+m^{2}(1-4\lambda \mu)(a_{\mu}a_{\lambda}-d_{\lambda}d_{\mu})\\
&[c_{\lambda},c_{\mu}]&=&0\\
&[d_{\lambda},d_{\mu}]&=
&m(1+2\lambda)c_{\mu}d_{\lambda}-
m(1+2\mu)d_{\mu}c_{\lambda}+m^{2}(1-4\lambda \mu)c_{\mu}c_{\lambda}
\end{array}
\]
These relations staisfy the $\lambda \leftrightarrow \mu$ exchange
symmetry. The generators are arranged in the $T$-matrices
\[
T_{\lambda}=\begin{pmatrix}
a_{\lambda} & b_{\lambda}\\
c_{\lambda} & d_{\lambda}
\end{pmatrix}, \qquad
T_{\mu}=\begin{pmatrix}
a_{\mu} & b_{\mu}\\
c_{\mu} & d_{\mu}
\end{pmatrix}
\]
and the associated coproduct and counit are
\[
\cop (T_{\lambda})=T_{\lambda}\ot T_{\lambda}, \qquad
\cnt (T_{\lambda})=\ident
\]
The quantum determinant
\[
D_{\lambda}=
a_{\lambda}d_{\lambda}-b_{\lambda}c_{\lambda}+ 
m(1-2\lambda)a_{\lambda}c_{\lambda} =
a_{\lambda}d_{\lambda}-c_{\lambda}b_{\lambda}-
m(1+2\lambda)c_{\lambda}d_{\lambda}
\]
satisfies $[D_{\lambda},D_{\mu}]\neq 0$ and $[D_{\lambda},T_{\lambda}]\neq
0$ unless $\lambda = 0$. The coalgebra structure is also invariant under
the $\lambda \leftrightarrow \mu$ exchange symmetry.

\end{document}